# KÖNIG'S INFINITY LEMMA FOR LOCALLY FINITE ULTRAMETRIC SPACES GENERATED BY LABELED TREES

Oleksiy Dovgoshey[1,2], Olga Rovenska[3]
[1]*Institute of Applied Mathematics and Mechanics of NASU, Slovyansk, Ukraine*
[2]*University of Turku, Turku, Finland*
[3]*Donbas State Engineering Academy, Kramatorsk, Ukraine*

***Abstract.*** We analyze the interplay between labeled trees and the ultrametric spaces they present. We provide characterizations of labeled trees that generate separable ultrametric spaces and those that generate locally finite ultrametric spaces. In particular, we establish an analog of König's Infinity Lemma for locally finite ultrametric spaces generated by labeled trees.

According to [11], any finite ultrametric space is isometrically describable by a Gurvich-Vyalyi representing tree. The Gurvich-Vyalyi representing trees form a subclass of finite trees equipped with a specific labeling of the vertex set. Trees with labeled vertices have been extensively studied, resulting in numerous contributions (see, e.g., the survey in [10]). The corresponding geometric interpretation of the Gurvich-Vyalyi representation [14] provides a framework for addressing various extremal problems related to finite ultrametric spaces [1-2,4]. Analogues of the Gurvich-Vyalyi representation and its geometric interpretation have recently been extended to the class of totally bounded ultrametric spaces [7].

Infinite trees with positive real-valued edge labelings are commonly referred to as *R*-trees. The structural relationships between finite subtrees of *R*-trees and finite monotone rooted trees are described in [3]. The categorical equivalence between trees and ultrametric spaces has been studied in [12, 13].

Ultrametric spaces generated by arbitrary nonnegative vertex labelings on finite and infinite trees were introduced in [3] and subsequently studied in [5, 6]. A characterization of totally bounded ultrametric spaces generated by labeled almost rays was provided in [9]. Furthermore, [8] contains a metric characterization of ultrametric spaces generated by labeled star graphs.

Our work [15] continues the aforementioned studies on the characterization of labeled trees generating ultrametric spaces, with a focus on the topological properties of these spaces. The aim is to describe the class of trees that admit labelings which generate separable ultrametric spaces. Additionally, we establish a characterization of such trees in terms of the local finiteness of the corresponding ultrametrics.

In what follows, we will use the terminology of [3].

By labeled tree $T(l)$ we mean a tree $T$ equipped with a labeling $l : V(T) \to \mathbb{R}^+$, where $V(T)$ is the vertex set of $T$. Let $T(l)$ be a labeled tree. We define a mapping $d_l$ : $V(T)\times V(T) \to \mathbb{R}^+$ by

$$d_l(u,v) = \begin{cases} 0, & \text{if } u = v \\ \max\limits_{w \in V(P)} l(w), & \text{if } u \neq v, \end{cases}$$

where $P$ is the path joining $u$ and $v$ in $T$. The mapping $d_l$ is an ultrametric on $V(T)$ iff the inequality $\max\{l(u), l(v)\}>0$ holds for all adjacent $u, v \in V(T)$.

The following two results were proved in [15].

**Theorem 1.** Let $T$ be a tree. Then the following statements are equivalent:

1. The vertex set of $T$ is countable.
2. All ultrametric spaces $(V(T), d_l)$ are separable.
3. There exists a labeling $l : V(T) \to \mathbb{R}^+$ such that $(V(T), d_l)$ is a separable ultrametric space.
4. For every star graph $S \subseteq T$ there exists a non-degenerate labeling $d_S : V(S) \to \mathbb{R}^+$ such that the ultrametric space $(V(S), d_{ls})$ is separable.

The König Infinity Lemma describes the finite connected graphs as graphs which have no vertices with infinite degree and contain no rays. The next theorem can be considered as an analog of this lemma for locally finite ultrametric spaces $(V(T), d_l)$.

**Theorem 2.** Let $T = T(l)$ be a separable ultrametric space generated by labeling $l : V(T) \to \mathbb{R}^+$. Then the following statements are equivalent:

1. The ultrametric space $(V(T), d_l)$ is locally finite.
2. The ultrametric spaces $(V(R), d_l|_{V(R)\times V(R)})$ and $(V(S), d_l|_{V(S)\times V(S)})$ are locally finite for all rays $R \subseteq T$ and all star graphs $S \subseteq T$.

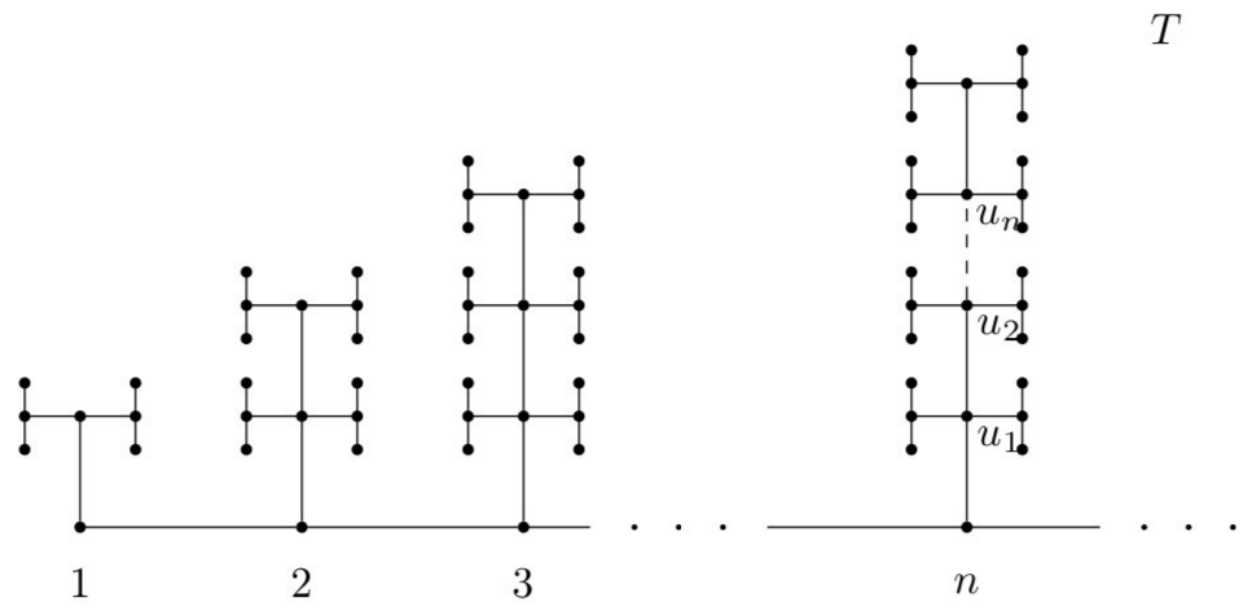

*Fig.1 – The tree $T$ does not contain any vertex of infinite degree and the symmetric difference of vertex sets $V(R_1)$, $V(R_2)$ is finite for any two rays $R_1 \subseteq T$, $R_2 \subseteq T$.*

**Example 1.** Let $T$ be the tree containing the ray $R = (v_1, v_2, \dots, v_n, \dots)$ with labeling $l_R : V(T) \to \mathbb{R}^+$ satisfying $l_R(v_n) = n$ for each $n \in \mathbb{N}$ (see Figure 1). If $l : V(T) \to \mathbb{R}^+$ is a non-degenerate labeling such that $l_R$ is a restriction of $l$ on the set $V(R)$, then the ultrametric space $(V(T), d_l)$ is locally finite by Theorem 2.

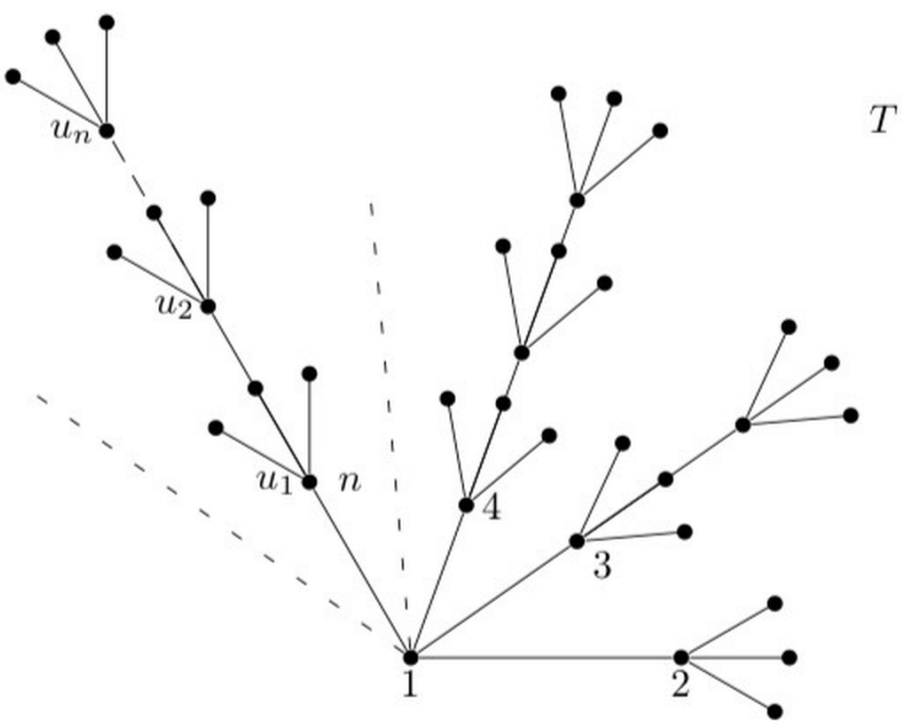

*Fig.2 – The tree T is rayless and contains exactly one vertex of infinite degree.*

**Example 2.** Let $T$ be the tree containing the vertex $v$ with $degT(c) = \aleph_0$, let $S$ be the star graph induced in $T$ by $c$ and all vertices adjacent to $c$, and let $l_S : V(S) \to \mathbb{R}^+$ be an injective labeling such that $l_S(V(S)) = \mathbb{N}$ (see Figure 2). If $l : V(T) \to \mathbb{R}^+$ is a non-degenerate labeling such that $l_S$ is a restriction of $l$ on the set $V(S)$, then the ultrametric space $(V(T), d_l)$ is locally finite by Theorem 2.